\documentclass[11,5pt, a4paper]{article}
\usepackage[utf8]{inputenc}
\usepackage{titlesec}
\usepackage[english]{babel}
\usepackage{graphicx}
\usepackage{float}
\usepackage{fancyhdr}
\usepackage[matrix,arrow,curve]{xy}
\usepackage{pstricks} 
\usepackage{amsmath,amsfonts,verbatim,afterpage,theorem,euscript,mathrsfs,amssymb}
\usepackage{amsfonts}
\usepackage{amssymb}
\usepackage{array}
\usepackage{dsfont}
\usepackage{hyperref}
\usepackage{authblk}
\newtheorem{Theorem}{Theorem}[section]
\newtheorem{Lemma}{Lemma}[section]

\newtheorem{Remark}{Remark}[section]
\newtheorem{Definition}{Definition}[section]

\numberwithin{equation}{section}

\def \vu{\textbf{u}}
\def \vv{\textbf{v}}

\begin{document}

\title{Micropolar fluids with initial angular velocities in non-homogeneous Sobolev spaces of order $-1/2$}
\author{Pedro Gabriel Fern\'andez-Dalgo \footnote{Basque Center for Applied Mathematics, BCAM, Bilbao, España \\ \textit{Email address:} \texttt{pfernandez@bcamath.org}} 
\footnote{Escuela de Ciencias Físicas y Matemáticas, Universidad de Las Américas, Vía a Nayón, C.P.170504, Quito, Ecuador}
 }
\date{}



\maketitle

\begin{abstract}
\noindent
In this paper, we investigate fractional energy methods for Micropolar fluids, starting with an initial angular velocity of negative Sobolev regularity. For the initial angular velocity assumption, we consider a non-homogeneous Sobolev norm of negative order. The regularity -1/2 studied here corresponds to the critical scaling of a simplified associated system, and the general framework can also be applied to the Boussinesq system with viscosity. Since our approach differs from those based on mild solutions and does not rely on a projected system, this work provides new tools for studying the Caffarelli-Kohn-Nirenberg theory of singularities in coupled variables within the Navier-Stokes equations.
\end{abstract}




\section{Introduction}
\noindent

Micropolar fluids are important examples of systems coupled to the Navier-Stokes equations, which model new physical phenomena thanks to the introduction of a new variable $\omega$ to analyze the effect of microrotations occurring in the fluid. The analysis of this coupled system leads to a better understanding of some experiments in which the Navier-Stokes equations are not sufficient to describe the fluid motion, particularly fluids containing polymeric additives. For a physical deduction and motivation of the micropolar model we refer to \cite{Eri,KMP83,CGOR}. Several studies about regularity of solutions for the micropolar system
\begin{equation*} \label{M}
(\text{M}) \left\{ \begin{array}{ll} \vspace{2mm} 
\partial_t \vu = \Delta \vu - (\vu \cdot \nabla) \vu  - \nabla p + \frac{1}{2} \nabla \wedge \omega ,  \\ \vspace{2mm}
\partial_{t} \omega = \Delta \omega - (\vu \cdot \nabla) \omega + \frac{1}{2} \nabla \wedge \vu - \omega +\nabla ( \nabla \cdot \omega )
, \\ \vspace{2mm}
\nabla \cdot  \vu=0
\end{array}
\right.     
\end{equation*}
suggest that the most restrictive hypotheses are on the velocity variable $\vu$ beacause of the presence of the term $(\vu \cdot \nabla) \vu$, which is nonlinear in $\vu$, and that the ideas applied to the Navier-Stokes equations can be extended to this system considering weaker hypothesis on $\omega$, we refer to \cite{JYuan, Bol, GR77, CLcr, CLps}. As in the case of the Navier-Stokes equations, the pressure $p$, introduced to save the incompressibility of $\vu$, is an auxiliary unknown which under very weak decaying assumptions on the whole space is overdetermined. In fact, taking the divergence in the equation for $\partial_t \vu$ and observing that $\nabla \cdot ( 
\nabla \wedge \omega) = 0$ we find $-\Delta p = \sum_{i,j} \partial_i \partial_j (u_i u_j)$.

The Micropolar equations (M) are not scaling invariant however the simplified system
\begin{equation*}
    (\text{sM}) \left\{ \begin{array}{ll} \vspace{2mm} 
    \partial_t \vu = \Delta \vu - (\vu \cdot \nabla) \vu  - \nabla p + \frac{1}{2} \nabla \wedge \omega ,  \\ \vspace{2mm}
    \partial_{t} \omega = \Delta \omega - (\vu \cdot \nabla) \omega
    \end{array}
    \right.     
    \end{equation*}
    is scaling invariant, and more precisely if $(\vu, \omega)$ is a solution of (sM) then $\vu _\lambda (x,t) = \lambda \vu(\lambda x, \lambda^2 t)$ and $\omega _\lambda (x,t) = \lambda^2 \omega(\lambda x, \lambda^2 t)$ form a solution of (sM). A critical space for the initial data $(\vu_0,\omega_0)$ is $\dot H ^{1/2} \times \dot H ^{-1/2}$.


For the mathematical study of the Micropolar fluids, Galdi and Rionero \cite{GR77} present the standard frameworks extending the most classical results for the Navier-Stokes equations, which means the study of weak solutions for $(\vu_0,\omega_0) \in L^2 \times L^2$ and mild solutions for $(\vu_0,\omega_0) \in L^3 \times L^3$ or $(\vu_0,\omega_0) \in \dot H^{1/2} \times H^{1/2}$, and an extension of the Leray's weak solutions framework is given in the paper of Boldrini and Rojas-Medar \cite{Bol} for a more general system including the magneto-hydrodynamics equations.
The current methods to construct strong or mild solutions of systems coupled to the Navier-Stokes equations use the fixed point theorem for contraction mappings in Banach spaces \cite{DeCu, GR77, JYuan}. We shall diverge of these methods by use energy balances with different order of regularity for each variable, in order to construct solutions starting from a well-behaved initial velocity and a microrotation data with non-homogeneous Sobolev regularity until -1/2. The computations begin from fractional energies, hence present work provides tools for the study of the CKN theory of singularities of the microrotation variable in a future work, the development of the CKN theory is due to Caffarelli, Kohn and Nirenberg,
and the founding works of Scheffer \cite{CKN, Sche76, Sche77}.

\section{Main results}
\noindent
We denote $R_i$ the i-th Riesz transform operator.
\begin{Theorem}
\label{mm} (Micropolar fluids)
Let $\vu_0 \in H^\tau(\mathbb{R}^3)$ be a divergence free initial velocity, and let $\omega_0 \in H^{\sigma}(\mathbb{R}^3)$ be an initial angular velocity.

\noindent
If $1/2 < \tau < 3/2$ and $\tau-1 < \sigma < 3/2$
then, there exists a positive time 
\begin{equation*}
    T_E = C_1 ( 1+ \| \vu_0 \|_{\dot{H}^\tau} + \| \omega_{0}  \|_{H^\sigma} )^{-\frac{2}{2\tau-1}} ,
\end{equation*}
for which
we find the existence of a solution $(\vu, \, \omega, \, p)$ on $[0,T_E]$ of the Cauchy problem (M) fulfilling

$(\textbf{P}1)$ $ \vu $ belongs to $ L^\infty ( (0, \, T_E), \, H^\tau ) \cap L^2 ( (0, \, T_E), \, \dot{H}^{\tau+1} ) $ and $\vu(t,\cdot) \rightarrow \vu_0$ when $t \to 0^+$ in $H^\tau$, 

$(\textbf{P}2)$ $ \omega $ belongs to $ L^\infty ( (0, \, T_E), \, H^{\sigma} ) \cap L^2 ( (0, \, T_E), \, \dot{H}^{\sigma+1} ) $  and $\omega(t,\cdot) \rightarrow \omega_0$ when $t \to 0^+$ in $H^\sigma$,
    
$(\textbf{P}3)$ $ p = \sum_{1 \leq i,j \leq 3} R_i R_j (u_i u_j)$.

\noindent
If $\tau=1/2$ and $-1/2 < \sigma < 3/2$, there exists $ \epsilon_{0} > 0 $ such that if $ \| \vu_0 \|_{H^{1/2}} + \| \omega_0 \|_{H^{\sigma}} < \epsilon_{0}$  then, there exists a positive time $T_E = T_E (\sigma, \tau )$, for which
we get the existence of a solution $(\vu, \, \omega, \, p)$ on $[0,T_E]$ of the Cauchy problem (M) fulfilling $(\textbf{P}1)$, $(\textbf{P}2)$ and $(\textbf{P}3)$.
\end{Theorem}

\noindent
In order to treat the case $\sigma = \tau-1 $, we consider an extra condition on the viscosity coefficients. Let us write,
\begin{equation*} \label{Mmunu}
(\text{M}^{\mu,\nu}) \left\{ \begin{array}{ll} \vspace{2mm} 
\partial_t \vu = \nu \Delta \vu - (\vu \cdot \nabla) \vu  - \nabla p + \frac{1}{2} \nabla \wedge \omega ,  \\ \vspace{2mm}
\partial_{t} \omega = \mu \Delta \omega - (\vu \cdot \nabla) \omega + \frac{1}{2} \nabla \wedge \vu - \omega +\nabla ( \nabla \cdot \omega )
, \\ \vspace{2mm}
\nabla \cdot  \vu=0, \\ \vspace{2mm}
\vu(0,\cdot)= \vu_{0} ( \cdot ), \phantom{spa}
\omega(0,\cdot)= \omega_{0} ( \cdot ).
\end{array}
\right.     
\end{equation*}

\noindent
In this case, we have the following result.
\begin{Theorem}
\label{mlc} (Micropolar fluids)
Let $\vu_0 \in H^\tau(\mathbb{R}^3)$ be a divergence free initial velocity, and let $\omega_0 \in H^{\sigma}(\mathbb{R}^3)$ be an initial angular velocity.

\noindent
There exists $N>0$, such that for all $\mu, \nu >0$ such that $ \mu \nu \geq N$,
if $1/2 \leq \tau < 3/2$ and $\sigma = \tau-1 $
then, there exists $ \epsilon_{0} > 0 $ such that if $ \| \vu_0 \|_{H^{1/2}} + \| \omega_0 \|_{H^{\sigma}} < \epsilon_{0}$  then, we can find a positive time $T_E = T_E (\sigma, \tau )$, for which we find a solution $(\vu, \, \omega, \, p)$ on $[0,T_E]$ of $(M^{ \mu, \nu })$ fulfilling
$(\textbf{P}1)$, $(\textbf{P}2)$ and $(\textbf{P}3)$.
\end{Theorem}

\noindent 
{\bf Related works.} Observing the method applied in \cite{BaChDa} (Theorem 5.6 of this book), we can expect uniqueness (or partial uniqueness at the endpoint regularity $\sigma=-1/2$) in the case where the initial angular velocity belongs to {\textit{homogeneous}} Sobolev spaces of negative order, however in our more general setting where the angular velocity belongs to {\textit{non-homogeneous}} Sobolev spaces of negative order, uniqueness is not immediately clear.
With respect to the regularity of these solutions we present the following remark.
\begin{Remark}
    By the properties described in $(\textbf{P}1)$ and $(\textbf{P}2)$ the product $\partial_t \vu \cdot \vu$ is well defined (see computations \eqref{rotterm1} and \eqref{rotterm2}) and in the sense of distributions we have
\begin{align*}
    \partial_t \, | \vu |^2 = \,  2 \vu \cdot \partial_t \vu
     = - 2 | \nabla \vu |^2 
    + \Delta(| \vu |^2 ) - \nabla \cdot \left( | \vu |^2 \vu  + 2 p \vu    \right) + (\nabla \wedge \omega_{\epsilon}) \cdot \vu .
\end{align*}
    However, if $\sigma<0$ then we are not able to demonstrate that the product $\partial_t \omega \cdot \omega$ is well defined, neither the existence of a non negative locally finite measure $\mu$ such that
\begin{align*}
    \mu =& - \partial_t \, | \omega|^2 - 2| \nabla \omega |^2 - 2 | \omega |^2 - 2| \nabla \cdot \omega |^2  \nonumber \\
    & + \Delta(| \omega |^2 ) - \nabla \cdot \left( | \omega |^2  \vu    \right)  + 2 \nabla ( (\nabla \cdot \omega) \omega ) + (\nabla \wedge \vu) \cdot \omega. \nonumber
\end{align*}    

\noindent
Thus, solutions appearing in Theorem \ref{mm} and Theorem \ref{mlc} belongs to a generalization of the class of partial suitable solutions introduced in Definition 1 in \cite{CLps}. By this reason, the analysis of regularity deserves a deep study where the Caffarelli, Kohn and Nirenberg appears naturally, it will be done in a future work.
\end{Remark}

\noindent
Recently, in the paper of Chamorro and Llerena \cite{CLps}, a notion of partial suitability is introduced in order to weaken hypothesis on the microrotations variable $\omega$ in the $\epsilon$-regularity result given in Theorem 1 in \cite{CLps}, where the authors assume in particular $\vu, \omega \in L^\infty_t L^2_x \cap L^2_t \dot H ^1_x $. However, in this paper we show the existence of solutions fulfilling only the weaker condition $\omega \in L^\infty_t H^{\sigma}_x \cap L^2_t \dot H ^{\sigma+1}_x $, where $-1/2 \leq \sigma < 0$, so that the result of Chamorro y Llerena must be studied in a more general functional setting for the microrotation variable.

\section{Proof of Theorems \ref{mm} and \ref{mlc}}

We begin by consider the classical approximated solutions obtained by mollification in the non-linear term, after the key idea is to mollify the time derivative of the approximated solutions and multiply this expression by a fractional Laplacian of the approximated solution.

\begin{Definition}
    We will write $D^s=(-\Delta)^{s/2}$ and $L^s = (I-\Delta)^{s/2}$, where $I$ is the identity operator, $\Delta$ is the Laplacian operator and $s \in \mathbb R$.
\end{Definition}

\subsection{The approximated micropolar system}
We fix the initial data $( \vu_0 , \omega_0) \in H^\tau(\mathbb{R}^3) \times H^\sigma(\mathbb{R}^3) $, where $\vu_0$ is divergence free.
Let us denote $\vv_{\epsilon} = \vu_{\epsilon}*\theta_{\epsilon}$.

\noindent
Consider $(\vu_{\epsilon}, \omega_{\epsilon},p_{\epsilon})$ the unique global solution of the mollified problem 
\begin{equation*}\label{Me1}
(\text{M}_{\epsilon}) \left\{ \begin{array}{ll}\vspace{2mm} 
\partial_t \vu_{\epsilon} = \Delta \vu_{\epsilon} - \vv_{\epsilon} \cdot \nabla \vu_{\epsilon}  - \nabla p_{\epsilon} + \frac{1}{2} \nabla \wedge \omega_{\epsilon} ,  \\ \vspace{2mm}
\partial_{t} \omega_{\epsilon} = \Delta \omega_{\epsilon} - \vv_{\epsilon} \cdot \nabla \omega_{\epsilon} + \frac{1}{2} \nabla \wedge \vu_{\epsilon} - \omega_{\epsilon} +\nabla ( \nabla \cdot \omega_{\epsilon} )
, \\ \vspace{2mm}
\nabla \cdot  \vu_{\epsilon} =0, \phantom{space}
\vu_{\epsilon}(0,\cdot)= \vu_{0} , \phantom{space}
\omega_{\epsilon}(0,\cdot)= \omega_{0} * \theta_\epsilon ,
\end{array}
\right.     
\end{equation*}
which satisfies $(\vu_{\epsilon},\omega_{\epsilon})$ belongs to $\mathcal{C}([0,+\infty), L^2(\mathbb{R}^d))\cap L^2((0,+\infty),\dot H^1(\mathbb{R}^d))$ and the pressure is given by the formula $ p_{\epsilon} = \sum_{1 \leq i, j \leq d} R_i R_j (v_{\epsilon, \lambda , i} u_{\epsilon, \lambda, j})$. This solution is smooth on $(0,+\infty) \times \mathbb R^3$. Moreover, the solution of this system satisfies the following classical energy balances:
\begin{align}
\label{enum}
    \partial_t \, | \vu_{\epsilon} |^2 = \,  2 \vu_{\epsilon} \cdot \partial_t \vu_{\epsilon}
     = - 2 | \nabla \vu_{\epsilon} |^2 
    + \Delta(| \vu_{\epsilon} |^2 ) - \nabla \cdot \left( | \vu_{\epsilon} |^2 \vv_{\epsilon}  + 2 p_{\epsilon} \vu_{\epsilon}    \right) + (\nabla \wedge \omega_{\epsilon}) \cdot \vu_{\epsilon}
\end{align}
and
\begin{align}
\label{enom}
    \partial_t \, | \omega_{\epsilon}|^2 = & \, 2 \omega_{\epsilon} \cdot \partial_t \omega_{\epsilon} \\
    =& - 2| \nabla \omega_{\epsilon} |^2 - 2 | \omega_{\epsilon} |^2 - 2| \nabla \cdot \omega_{\epsilon} |^2  \nonumber \\
    & + \Delta(| \omega_{\epsilon} |^2 ) - \nabla \cdot \left( | \omega_{\epsilon} |^2  \vv_{\epsilon}    \right)  + 2 \nabla ( (\nabla \cdot \omega_{\epsilon}) \omega_{\epsilon} ) + (\nabla \wedge \vu_{\epsilon}) \cdot \omega_{\epsilon}. \nonumber
\end{align} 
\noindent
We can integrate in space and time on an arbitrary interval $(0,t)$ the equations \eqref{enum} and \eqref{enom} in order to obtain
\begin{align}
\label{enui1}
    \| \vu_{\epsilon} (t)\|^2_{L^2} 
    + 2 \int_0^t \int | \nabla \vu_{\epsilon} |^2 = \| \vu_{\epsilon}(0)\|^2_{L^2}  + \int_0^t \int (\nabla \wedge \omega_{\epsilon}) \cdot \vu_{\epsilon}
\end{align}
and
\begin{align}
\label{enoi1}
    \| \omega_{\epsilon} (t) \|^2_{L^2} 
    + 2 \int_0^t \int | \nabla \omega_{\epsilon} |^2 + | \omega_{\epsilon} |^2 + | \nabla \cdot \omega_{\epsilon} |^2 \,  
     = \| \omega_{\epsilon} (0) \|^2_{L^2} + \int_0^t \int (\nabla \wedge \vu_{\epsilon}) \cdot \omega_{\epsilon}.
\end{align}   

\noindent
Observe that when $\sigma<0$, the equation \eqref{enoi1} is not useful in order to take the limit when $\epsilon$ goes to $0$ as $\omega_0 \in H^\sigma$ does not implies $\omega_0 \in L^2$.

We shall consider energy balances where a non local operator intervenes in order to control fractional derivatives of the velocity following the parameter $\tau$,
\begin{align}
\label{lapue}
    \left( D^\tau \vu_{\epsilon} * \theta_\kappa \right) \cdot  \left( \partial_t D^\tau \vu_{\epsilon} * \theta_\kappa \right) =& \, \, \, \, \, \, \left( D^\tau \vu_{\epsilon} * \theta_\kappa \right) \cdot  \Delta \left( D^\tau \vu_{\epsilon} * \theta_\kappa \right) \\
    &- \left( D^\tau \vu_{\epsilon} * \theta_\kappa \right) \cdot   \left( D^\tau \nabla \cdot ( \vv_{\epsilon} \otimes \vu_{\epsilon} ) * \theta_\kappa \right) \nonumber  \\
    &- \left(  D^\tau \vu_{\epsilon} * \theta_\kappa \right) \cdot  D^\tau \nabla  R_i R_j (v_{\epsilon, i} u_{\epsilon, j} ) * \theta_\kappa  \nonumber \\
    &+ \left( D^\tau \vu_{\epsilon} * \theta_\kappa \right) \cdot  \frac{1}{2} D^\tau \left(  \nabla \wedge \omega_{\epsilon} \right) * \theta_\kappa ,  \nonumber
\end{align} 
and fractional derivatives of the angular velocity following the parameter $\sigma$,
\begin{align}
\label{lapoe}
    \left( L^{\sigma} \omega_{\epsilon} * \theta_\kappa \right) \cdot  \left( \partial_t L^{\sigma} \omega_{\epsilon} * \theta_\kappa \right) &=   \, \, \, \, \, \, \left( L^{\sigma} \omega_{\epsilon} * \theta_\kappa \right) \cdot  \Delta \left( L^{\sigma} \omega_{\epsilon} * \theta_\kappa \right) \\
    &-  \left( L^{\sigma} \omega_{\epsilon} * \theta_\kappa \right) \cdot   \left( L^{\sigma} \nabla \cdot ( \vv_{\epsilon} \otimes \omega_{\epsilon} ) * \theta_\kappa \right) \nonumber  \\
    &+ \left( L^{\sigma} \omega_{\epsilon} * \theta_\kappa \right) \cdot \frac{1}{2} L^{\sigma} \left(  \nabla \wedge \vu_{\epsilon} - 2 \omega_{\epsilon} + 2 \nabla ( \nabla \cdot \omega_{\epsilon} ) \right) *\theta_\kappa  . \nonumber
\end{align}  

\noindent
The convolution by $\theta_\kappa$ is introduced in \eqref{lapoe} in order to have $L^{\sigma} \omega_\epsilon * \theta_\kappa \in L^2 H^{1}$ and $ \partial_t L^\sigma \omega_\epsilon * \theta_\kappa \in L^2 H^{-1} $. Moreover, we have $D^{\tau} \vu_\epsilon * \theta_\kappa \in L^2 H^{1}$ and $ \partial_t D^\tau \vu_\epsilon * \theta_\kappa \in L^2 H^{-1} $. Then,
integration of \eqref{lapue} on the whole space and over an arbitrary interval $(0,t)$ gives
\begin{align}
\label{intue}
    \|  D & ^\tau \left( \vu_{\epsilon} * \theta_\kappa \right) (t) \|_{L^2}^2 + 2 \int_0^t \int  | \nabla D^\tau \left( \vu_{\epsilon} * \theta_\kappa \right) |^2  \\
    \leq& \, \|  D^\tau  \left( \vu_{\epsilon} * \theta_\kappa \right) (0) \|_{\dot{H}^\tau}^2 + \int_0^t \int \Delta \, | D^\tau \left( \vu_{\epsilon} * \theta_\kappa \right) |^2  \nonumber \\
    & -2 \int_0^t \int \left( D^\tau \nabla \cdot ( \vv_{\epsilon} \otimes \vu_{\epsilon} ) * \theta_\kappa \right) \cdot \left( D^\tau \vu_{\epsilon} * \theta_\kappa \right)  \\
    & -2\int_0^t \int D^\tau \nabla \left( R_i R_j ( v_{\epsilon, i} u_{\epsilon, j} ) * \theta_\kappa  \right) \cdot \left( D^\tau \vu_{\epsilon} * \theta_\kappa \right)  \nonumber \\
    & + \int_0^t \int \left( D^\tau \nabla \wedge \omega_{\epsilon} * \theta_\kappa \right) \cdot \left( D^\tau \vu_{\epsilon} * \theta_\kappa \right) \nonumber
\end{align} 
and $\int_0^t \int \Delta \, | D^\tau \left( \vu_{\epsilon} * \theta_\kappa \right) |^2 = 0 $. Moreover, from \eqref{lapoe} we get
\begin{align}
\label{intoe}
    &\, \, \, \| L^{\sigma}  ( \omega_{\epsilon}  *  \theta_\kappa ) (t) \|_{L^2}^2 \\
    &+2 \int_0^t \int | \nabla L^{\sigma} \left( \omega_{\epsilon} * \theta_\kappa \right) |^2 \nonumber
    + | L^{\sigma} \left( \omega_{\epsilon} * \theta_\kappa \right) |^2 + | \nabla \cdot L^{\sigma} \left( \omega_{\epsilon} * \theta_\kappa \right) |^2 \, \, dx ds  \nonumber  \\
    & \, \, \, \, \, \, \, \leq \|  L^{\sigma}  \left( \omega_{\epsilon} * \theta_\kappa \right) (0) \|_{\dot{H}^\sigma}^2 + \int_0^t \int \Delta \, |L^{\sigma} \left( \omega_{\epsilon} * \theta_\kappa \right) |^2  \nonumber \\
    & \, \, \, \, \, \, \, -\int_0^t \int \left( L^{\sigma} \nabla \cdot ( \vv_{\epsilon} \otimes \omega_{\epsilon} ) * \theta_\kappa \right) \cdot \left( L^{\sigma} \omega_{\epsilon} * \theta_\kappa \right)  \nonumber \\
    & \, \, \, \, \, \, \, + \int_0^t \int \left( L^{\sigma} \nabla \wedge \vu_{\epsilon} * \theta_\kappa \right) \cdot \left( L^{\sigma} \omega_{\epsilon} * \theta_\kappa \right) , \nonumber
\end{align} 
and $ \int_0^t \int \Delta \, |L^{\sigma} \left( \omega_{\epsilon} * \theta_\kappa \right) |^2 = 0 $.

\subsection{Controlling the velocity in $H^\tau$}
We search to estimate the right hand side of \eqref{intue} absorbing the less regular quantities.

We begin by bound the $L^\infty_t L^2_x$ norm of $\vu$. For that, from \eqref{enui1} we observe that the term to be controlled comes from the coupled rotational. \noindent
We analyze now the coupled rotational. First, suppose $ -1/2 \leq \sigma \leq 1$ and take $ s \in [\tau-1/2, \tau] $ such that
$-\tau+s  \leq \sigma \leq 1 - \tau +s $.
Then,
\begin{align*}
    \int_0^t \int (\nabla \wedge \omega_{\epsilon}) \cdot \vu_{\epsilon} \leq \int_0^t \| \omega \|_{\dot{H}^{1-\tau+s}} \| \vu \|_{\dot{H}^{\tau-s}}
\end{align*}
and $ \max \{ 1/2, \sigma \} \leq 1- \tau+s \leq \sigma+1 $. Thus, using interpolation we find
\begin{align*}
    \| \omega \|_{\dot{H}^{1 -\tau +s}} \| \vu \|_{\dot{H}^{\tau-s}}
    \leq & \| \omega \|_{H^{1 -\tau +s} } \| \vu \|_{H^{\tau}} \nonumber \\
    & \| \omega \|_{H^{\sigma}}^{ \tau + \sigma -s } \| \omega \|_{H^{\sigma+1}}^{ 1 -\tau - \sigma +s } \| \vu \|_{H^{\tau}} ,\nonumber
\end{align*}
from where it can be obtained that
\begin{align}
\label{rotterm1}
    \int_0^t \int (\nabla \wedge \omega_{\epsilon}) \cdot \vu_{\epsilon} \leq C_\delta \int_0^t \| \vu_{\epsilon} \|_{H^{\tau}}^{ 2 } + \| \omega_{\epsilon} \|_{H^{\sigma}}^{ 2 } \, \, ds   + \delta \int_0^t \| \omega_{\epsilon} \|_{H^{\sigma+1}}^2 .
\end{align}
Now, suppose $1 < \sigma$, then
\begin{align*}
    \int_0^t \int (\nabla \wedge \omega_{\epsilon}) \cdot \vu_{\epsilon} \leq \int_0^t \| \omega_\epsilon \|_{\dot{H}^{1}} \| \vu_\epsilon \|_{L^2} \leq \int_0^t \| \omega_\epsilon \|_{H^{\sigma}} \| \vu_\epsilon \|_{L^2}
\end{align*}
so that
\begin{align}
\label{rotterm2}
    \int_0^t \int (\nabla \wedge \omega_{\epsilon}) \cdot \vu_{\epsilon} \leq C \int_0^t \| \vu_{\epsilon} \|_{L^2}^{ 2 } + \| \omega_{\epsilon} \|_{H^{\sigma}}^{ 2 } \, \, ds  .
\end{align}

The following step is to control $\vu$ in the $L^\infty_t \dot{H}^\tau$ seminorm.

\noindent
We analyze first the pressure part. As $ \tau < 3/2 $, there exists $ 0 < s < 1/2 $ such that $ \tau < 3/2 - s $,
\begin{align*}
    -\int_0^t \int D^{\tau} \nabla ( R_i R_j ( v_{\epsilon , i} u_{\epsilon, j}) * \theta_\kappa ) \cdot \left( D^{\tau} \vu_{\epsilon} * \theta_\kappa \right) \leq \int_0^t \| \vv_\epsilon \otimes \vu_\epsilon \|_{\dot{H}^{\tau}} \| \vu_\epsilon \|_{\dot{H}^{\tau+1}} . 
\end{align*}
By the product laws, as $ 0 < 3/2 - s < 3/2 $ and $ 0 \leq \tau + s $, we obtain
\begin{align*}
    \| \vv_{\epsilon} \otimes \vu_\epsilon \|_{\dot{H}^{\tau}} \leq \| \vv_\epsilon \|_{\dot{H}^{ 3/2 - s }} \| \vu_\epsilon \|_{\dot{H}^{ \tau + s }} + \| \vu_\epsilon \|_{\dot{H}^{ 3/2 - s }} \| \vv_\epsilon \|_{\dot{H}^{ \tau + s }} \leq \| \vu_\epsilon \|_{\dot{H}^{ 3/2 - s }} \| \vu_\epsilon \|_{\dot{H}^{ \tau + s }} .
\end{align*}
We have $ \tau \leq 3/2-s \leq \tau +1 $, or equivalently $ 1/2 - s \leq \tau \leq 3/2 - s $, thus we can interpolate to write
\begin{align*}
    \| \vu_\epsilon \|_{\dot{H}^{ 3/2 - s }} \| \vu_\epsilon \|_{\dot{H}^{ \tau + s }}   &\leq 
    \| \vu_\epsilon
    \|_{\dot{H}^{\tau} }^{ s + \tau -1/2 }
    \| \vu_\epsilon \|_{\dot{H}^{ \tau+1 }}^{ - s -\tau + 3/2 } \| \vu_\epsilon
    \|_{\dot{H}^{\tau} }^{ 1 - s }
    \| \vu_\epsilon \|_{\dot{H}^{ \tau+1 }}^{ s } \\
    &= \| \vu_\epsilon
    \|_{\dot{H}^{\tau} }^{ \tau +1/2 }
    \| \vu_\epsilon \|_{\dot{H}^{ \tau+1 }}^{ -\tau + 3/2 } ,
\end{align*}
hence, one gets
\begin{align*}
    -\int_0^t \int D^{\tau} \nabla ( R_i R_j ( v_{\epsilon , i} u_{\epsilon, j}) * \theta_\kappa ) \cdot \left( D^{\tau} \vu_{\epsilon} * \theta_\kappa \right) \leq \int_0^t \| \vu_\epsilon \|_{\dot{H}^{\tau}}^{ \frac{2\tau +1}{2} }  \| \vu_\epsilon \|_{\dot{H}^{\tau+1}}^{\frac{5-2\tau}{2}} .
\end{align*}
Then, as $1/2 < \tau$ we find
\begin{align}
    -&\int_0^t \int D^{\tau} \nabla ( R_i R_j ( v_{\epsilon , i} u_{\epsilon, j}) * \theta_\kappa ) \cdot \left( D^{\tau} \vu_{\epsilon} * \theta_\kappa \right) \nonumber \\ &\leq C_\delta \int_0^t \| \vu_\epsilon \|_{\dot{H}^{\tau}}^{2\left( \frac{2\tau +1}{2\tau -1 } \right) }  + \delta \int_0^t \| \vu_\epsilon \|_{\dot{H}^{\tau+1}}^2 .
\end{align}
For $\tau=1/2$, we obtain
\begin{align}
\label{pree}
    -\int_0^t \int D^{\tau} \nabla ( R_i R_j ( v_{\epsilon , i} u_{\epsilon, j}) * \theta_\kappa ) \cdot \left( D^{\tau} \vu_{\epsilon} * \theta_\kappa \right) \leq C \int_0^t \| \vu_\epsilon \|_{\dot{H}^{\tau}}^{1}  \| \vu_\epsilon \|_{\dot{H}^{\tau+1}}^{2} .
\end{align}

Now, we analyze the coupling part. We consider two cases, $\tau - 1 < \sigma < 2\tau $ and $ \tau < \sigma \leq 2\tau +1 $. 

When $\tau - 1 < \sigma \leq 2\tau $ we write
\begin{align*}
    \int_0^t \int (D^{\tau} \nabla \wedge \omega_{\epsilon}) \cdot \left( D^{\tau} \vu_{\epsilon} * \theta_\kappa \right) \leq \int_0^t \| \omega \|_{\dot{H}^{\sigma+1}} \| \vu \|_{\dot{H}^{2\tau-\sigma}}.
\end{align*}
Thus, using the fact that $ 0 \leq 2 \tau - \sigma \leq \tau+1 $, we find by interpolation
\begin{align*}
    \| \omega \|_{\dot{H}^{\sigma+1}} \| \vu \|_{\dot{H}^{2\tau-\sigma}} \leq \| \omega \|_{\dot{H}^{\sigma+1}} \| \vu \|_{L^2}^{1-\frac{2\tau-\sigma}{\tau+1}} \| \vu \|_{\dot{H}^{\tau+1}}^{\frac{2\tau-\sigma}{\tau+1}} .
\end{align*}
Then, as $ \tau-1 < \sigma $ we obtain
\begin{align}
        \int_0^t \int ( D^{\tau} \nabla \wedge \omega_{\epsilon}) \cdot \left( D^{\tau} \vu_{\epsilon} * \theta_\kappa \right) \leq C_\delta \int_0^t \| \vu_{\epsilon} \|_{L    ^2}^{ 2 } + \delta \int_0^t \| \omega_{\epsilon} \|_{\dot{H}^{\sigma+1}}^{ 2 } + \| \vu_{\epsilon} \|_{\dot{H}^{\tau+1}}^2 \, \, ds  .
\end{align}

\noindent
\begin{Remark}
    If $ \sigma = \tau -1 $, we write for $ \mu, \nu >0 $,
    \begin{align}
    \label{absorvenu}
        \int_0^t \int ( D^{\tau} \nabla \wedge \omega_{\epsilon}) \cdot \left( D^{\tau} \vu_{\epsilon} * \theta_\kappa \right) \leq &  C \int_0^t ( \mu^{\frac{1}{2}} \| \omega_{\epsilon} \|_{\dot{H}^{\sigma+1}} ) ( \nu^{\frac{1}{2}} \| \vu_{\epsilon} \|_{\dot{H}^{\tau+1}} ) .
\end{align}
\end{Remark}

In the second case, $ \tau < \sigma \leq 2\tau +1 $, we write
\begin{align*}
    \int_0^t \int ( D^{\tau} \nabla \wedge \omega_{\epsilon}) \cdot \left( D^{\tau} \vu_{\epsilon} * \theta_\kappa \right) \leq \int_0^t \| \omega \|_{\dot{H}^{\sigma}} \| \vu \|_{\dot{H}^{2\tau-\sigma+1}}.
\end{align*}
As we have $ 0 \leq 2\tau -\sigma +1 \leq \tau+1 $, using interpolation we find
\begin{align*}
    \| \omega \|_{\dot{H}^{\sigma}} \| \vu \|_{\dot{H}^{2\tau-\sigma+1}} \leq \| \omega \|_{\dot{H}^{\sigma}} \| \vu \|_{L^2}^{1-\frac{2\tau-\sigma+1}{\tau+1}} \| \vu \|_{\dot{H}^{\tau+1}}^{\frac{2\tau-\sigma+1}{\tau+1}} .
\end{align*}
Then, using the fact that $ \tau < \sigma $,
\begin{align}
        \int_0^t \int ( D^{\tau} \nabla \wedge \omega_{\epsilon}) \cdot \left( D^{\tau} \vu_{\epsilon} * \theta_\kappa \right)  \leq C_\delta \int_0^t \| \vu_{\epsilon} \|_{L    ^2}^{ 2 } + \| \omega_{\epsilon} \|_{H^{\sigma}}^{ 2 } \, \, ds  + \delta \int_0^t \| \vu_{\epsilon} \|_{\dot{H}^{\tau+1}}^2.
\end{align}

\subsection{Controlling the angular velocity in $H^\sigma$}

\noindent
We consider first terms arising from the non-linear part,
\begin{align*}
    -\int_0^t \int ( L^{\sigma} \nabla \cdot  ( \vv \otimes \omega) * \theta_\kappa ) \cdot \left( L^{\sigma} \omega_{\epsilon} *  \theta_\kappa \right) \leq \int_0^t \| \vv_{\epsilon} \otimes \omega_{\epsilon} \|_{H^{\sigma}} \| \omega_{\epsilon} \|_{H^{\sigma+1}}. 
\end{align*}
\noindent
As we have $1/2 \leq \tau \leq 3/2 $ and $-1/2 \leq \sigma \leq 3/2 $, then
\begin{equation}
\label{ob1}
    \tau - \left( \frac{3-2\tau}{2} \right) \leq -1/2 \leq \sigma \leq 3/2 \leq \tau +  \left( \frac{2\tau+1}{2} \right) .
\end{equation}
Observe that in view of \eqref{ob1} we have 
\begin{equation*}
\label{constraint}
   0 < \frac{ 3 + 2\sigma }{4} < 3/2 \phantom{spa} \text{and} \phantom{spa}  \max \left \{ \sigma, \tau \right \}  \leq \frac{ 3 + 2\sigma }{4} \leq \min \left \{ \sigma +1, \tau+1 \right \} .  
\end{equation*}
Then, by the product laws (see \cite{LR24}) and by interpolation we can obtain
\begin{align}
\label{pre}
    \| \vv_{\epsilon} \otimes \omega_{\epsilon} \|_{H^{ \sigma }} \leq & \| \vu_{\epsilon} \|_{H^{\frac{ 3 + 2\sigma }{4} }} \| \omega_{\epsilon} \|_{\dot{H}^{\frac{ 3 + 2\sigma }{4}}} + \| \vu_{\epsilon} \|_{\dot{H}^{\frac{ 3 + 2\sigma }{4} }} \| \omega_{\epsilon} \|_{H^{\frac{ 3 + 2\sigma }{4}}} \\
    \leq &
    \left( \| \vu \|_{L^2} + \| \vu_{\epsilon}
    \|_{\dot{H}^{\tau} }^{ \tau + \frac{ 1 -2\sigma  }{4} }
    \| \vu_{\epsilon} \|_{\dot{H}^{\tau+1}}^{ \frac{ 3 +2\sigma }{4} - \tau } \right) \left( \| \omega_{\epsilon}
    \|_{H^{\sigma} }^{ \frac{2\sigma + 1 }{4} }
    \| \omega_{\epsilon} \|_{H^{\sigma+1}}^{ \frac{ 3 - 2 \sigma }{4} } \right) . \nonumber
\end{align}
Moreover, using $1/2 < \tau$ one gets for $\delta>0$,
\begin{align*}
    \| \vu_{\epsilon}
    \|_{\dot{H}^{\tau} }^{ \tau + \frac{ 1 -2\sigma }{4} }
    \| \vu_{\epsilon} \|_{\dot{H}^{\tau+1}}^{ \frac{ 3 + 2\sigma  }{4} -\tau } \| \omega_{\epsilon}
    \|_{H^{\sigma} }^{ \frac{ 2\sigma +1 }{4}}
    \| \omega_{\epsilon} \|_{H^{\sigma+1}}^{ \frac{ 7 - 2\sigma }{4} } &\leq C_\delta \left( \| \vu_{\epsilon} \|_{\dot{H}^{\tau}}^{ 2\left( \frac{ 2\tau +1 }{2\tau -1} \right) } + \| \omega_{\epsilon} \|_{H^{\sigma}}^{ 2\left( \frac{2\tau +1 }{2\tau -1} \right) } \right)  \nonumber \\
    &\,\,\, + \delta \left( \| \vu_{\epsilon} \|_{\dot{H}^{\tau+1}}^2 + \| \omega_{\epsilon} \|_{H^{\sigma+1}}^2 \right) .
\end{align*}
Thus, we have found
\begin{align}
\label{subc}
    -\int_0^t \int ( L^{\sigma} \nabla & \cdot  ( \vv \otimes \omega) * \theta_\kappa ) \cdot \left( L^{\sigma} \omega_{\epsilon} *  \theta_\kappa \right) \nonumber\\
    &\leq C_\delta \int_0^t \| \vu \|_{L^2}^2 + \| \omega \|_{H^\sigma}^2 + \| \vu_{\epsilon} \|_{\dot{H}^{\tau}}^{ 2\left( \frac{ 2\tau +1 }{2\tau -1} \right) } + \| \omega_{\epsilon} \|_{H^{\sigma}}^{ 2\left( \frac{2\tau +1 }{2\tau -1} \right) } \, \, ds \nonumber\\
    &\,\,\, + \delta \int_0^t \| \vu_{\epsilon} \|_{\dot{H}^{\tau+1}}^2 + \| \omega_{\epsilon} \|_{H^{\sigma+1}}^2  \, \, ds . 
\end{align}
\noindent
If $\tau= 1/2$ then
\begin{align*}
    \| \vv_{\epsilon} \otimes \omega_{\epsilon} \|_{H^{ \sigma }} \leq & \| \vu_{\epsilon} \|_{H^{\frac{ 3 + 2\sigma }{4} }} \| \omega_{\epsilon} \|_{\dot{H}^{\frac{ 3 + 2\sigma }{4}}} + \| \vu_{\epsilon} \|_{\dot{H}^{\frac{ 3 + 2\sigma }{4} }} \| \omega_{\epsilon} \|_{H^{\frac{ 3 + 2\sigma }{4}}} \\
    \leq &
    \left( \| \vu \|_{L^2} + \| \vu_{\epsilon}
    \|_{\dot{H}^{\tau} }^{ \frac{ 3 -2\sigma  }{4} }
    \| \vu_{\epsilon} \|_{\dot{H}^{\tau+1}}^{ \frac{ 1 +2\sigma }{4} } \right) \left( \| \omega_{\epsilon}
    \|_{H^{\sigma} }^{ \frac{2\sigma + 1 }{4} }
    \| \omega_{\epsilon} \|_{H^{\sigma+1}}^{ \frac{ 3 - 2 \sigma }{4} } \right) \nonumber
\end{align*}
so that
\begin{align}
\label{post1/2}
    -&\int_0^t \int ( L^{\sigma} \nabla \cdot  ( \vv \otimes \omega) * \theta_\kappa ) \cdot \left( L^{\sigma} \omega_{\epsilon} *  \theta_\kappa \right) \\
    & \leq C \int_0^t \left( \| \vu_\epsilon \|_{H^ {\tau} } + \| \omega_\epsilon \|_{H^ {\sigma} } \right) \left( \| \vu_\epsilon \|_{ \dot{H}^ {\tau+1} }^{ 2 } + \| \omega_\epsilon \|_{H^ {\sigma} }^{ 2 } + \| \omega_\epsilon \|_{H^ {\sigma+1} }^{ 2 } \right) .
\end{align}
In the same way we have obtained \eqref{rotterm1}, we can get
\begin{align}
\label{rotterm3}
    \int_0^t \int (\nabla \wedge \vu_{\epsilon}) \cdot \omega_{\epsilon} \leq C_\delta \int_0^t \| \vu_{\epsilon} \|_{H^{\tau}}^{ 2 } + \| \omega_{\epsilon} \|_{H^{\sigma}}^{ 2 } \, \, ds  + \delta \int_0^t \| \omega_{\epsilon} \|_{H^{\sigma+1}}^2 
\end{align}
and similarly to \eqref{rotterm2}, we find
\begin{align}
\label{rotterm4}
    \int_0^t \int (\nabla \wedge \vu_{\epsilon}) \cdot \omega_{\epsilon} \leq C \int_0^t \| \vu_{\epsilon} \|_{L^2}^{ 2 } + \| \omega_{\epsilon} \|_{H^{\sigma}}^{ 2 } \, \, ds .
\end{align}

\noindent
Now, observe that the left hand side of \eqref{intoe} permit to control the $L^2 H^{\sigma+1}$ as
\begin{equation*}
    \int_0^t \int | \nabla L^{\sigma}  \omega_{\epsilon}  |^2 + | L^{\sigma} \omega_{\epsilon}  |^2 \geq c_3 \int_0^t \| \omega_\epsilon \|_{H^{\sigma+1}}^2.
\end{equation*}
Thus, we let $\kappa$ goes to $+\infty$ and we take $\delta$ small enough in order to obtain:

\noindent
Conclusion 1 : if $1/2 < \tau < 3/2$ and $ \tau-1 < \sigma < 3/2$ (or $ \mu \nu $ large enough and $\tau-1 = \sigma$), we get from \eqref{intue}-\eqref{rotterm4},
\begin{align}
\label{fin}
    \| \vu & (t)\|^2_{H^\tau} + \| \omega_\epsilon(t) \|_{H^\sigma}^2 + c \int_0^t \| \vu_\epsilon \|_{\dot{H}^{\tau+1}}^2 + \| \omega_\epsilon \|_{H^{\sigma+1}}^2 \, \, ds  \\
    \leq & \, \| \vu_{0} \|^2_{H^\tau} + \| \omega_{0, \epsilon} \|_{H^{\sigma}}^2 + C \int_0^t \| \vu_\epsilon \|_{H^ {\tau} }^2 + \| \omega_\epsilon \|_{H^ {\sigma} }^2 + \| \vu_\epsilon \|_{ H^ {\tau} }^{ 2\left( \frac{ 2\tau +1 }{2\tau -1} \right) } + \| \omega_\epsilon \|_{H^ {\sigma} }^{ 2\left( \frac{ 2\tau +1 }{2\tau -1} \right) } \, \, ds  . \nonumber
\end{align} 
\noindent
Conclusion 2 : if $\tau=1/2$ and $1/2 \leq \sigma < 3/2$ from \eqref{pree} and \eqref{post1/2},
\begin{align*}
    \| \vu_\epsilon & (t)\|^2_{H^\tau} + \| \omega_\epsilon(t) \|_{H^\sigma}^2 + c \int_0^t \| \vu_\epsilon \|_{\dot{H}^{\tau+1}}^2 + \| \omega_\epsilon \|_{H^{\sigma+1}}^2 \, \, ds  \\
    \leq & \, \| \vu_{0} \|^2_{H^\tau} + \| \omega_{0, \epsilon} \|_{H^{\sigma}}^2 + C \int_0^t \| \vu_{\epsilon} \|_{H^{\tau}}^2 + \| \omega_{\epsilon} \|_{H^{\sigma}}^2 \, \, ds  \nonumber \\
    & + C \int_0^t ( \| \vu_{\epsilon} \|_{H^{\tau}} + \| \omega_{\epsilon} \|_{H^{\sigma}} ) \left( \| \vu_{\epsilon} \|_{\dot{H}^{\tau+1}}^2 + \| \omega_{\epsilon} \|_{ H^{\sigma} }^2 + \| \omega_{\epsilon} \|_{ H^{\sigma+1} }^2 \right) \, \, ds .  \nonumber
\end{align*}
Hence, when $\tau=1/2$ and $ 1/2 \leq \sigma < 3/2$, under the assumption $ \| \vu_\epsilon (s, \cdot) \|_{H^{\tau}} + \| \omega_\epsilon (s, \cdot) \|_{H^{\sigma}} < \epsilon_0 $ with $C_2 \, \epsilon_0 <\frac 1 8$ where $C_2>0$ is a fixed constant, we get
\begin{align}
\label{fin2}
    \| \vu & (t)\|^2_{H^\tau} + \| \omega_\epsilon(t) \|_{H^\sigma}^2 + c \int_0^t \| \vu_\epsilon \|_{\dot{H}^{\tau+1}}^2 + \| \omega_\epsilon \|_{H^{\sigma+1}}^2 \, \, ds \\
    \leq & \, \| \vu_{0} \|^2_{H^\tau} + \| \omega_{0, \epsilon} \|_{H^{\sigma}}^2 + C \int_0^t \| \vu_{\epsilon} \|_{H^{\tau}}^2 + \| \omega_{\epsilon} \|_{H^{\sigma}}^2 \, \, ds . \nonumber 
\end{align}

\subsection{Passage to the limit of the approximated micropolar fluids}

\noindent
We use the following Gr\"onwall Lemma to obtain uniform estimates (we refer to Lemma 3.5 in \cite{PF} or \cite{PF_PG_arma}).

\begin{Lemma}
\label{lem5} 
Consider a continuous non-negative function $\alpha$ defined on $[0,T)$ which satisfies, for $A,B \in (0, + \infty)$  and $ b \in [1, \infty) $,
$$ \alpha(t)\leq A + B\int_0^t \alpha(s) + \alpha(s)^b\, ds.$$

\noindent
Then,
\begin{itemize}
    \item if $b>1$, we let $T_1 \in (0,T)$ and $T_0=\min \left\{ T_1, \frac 1{3^b B (A ^{b-1}+ (BT_1)^{b-1})} \right\} $. Then, we have, for every $t\in [0,T_0]$,  $\alpha(t)\leq  3A$.
    
    \item if $b=1$ we have for every $t\in [0,\frac{1}{4B}]$, the estimate  $\alpha(t)\leq  2A$.
\end{itemize}

\end{Lemma}

\noindent
Applying Lemma \ref{lem5} to the inequalities \eqref{fin} with $b=\frac{2\tau+1}{2\tau-1}$, we find that there exists a constant $ c_1 \geq 1$ such that if $T_0$ satisfies
\begin{equation*}
    \label{t1}
    c_1 ( 1+ \| \vu_0 \|_{\dot{H}^\tau} + \| \omega_{0}  \|_{H^\sigma} )^{\frac{2}{2\tau-1}} T_0 < 1 ,
\end{equation*}
then
\begin{align}
\label{unt0}
    \| \vu_\epsilon & (t)\|^2_{H^\tau} + \| \omega_\epsilon(t) \|_{H^\sigma}^2 + \int_0^t \| \vu_\epsilon \|_{\dot{H}^{\tau+1}}^2 + \| \omega_\epsilon \|_{H^{\sigma+1}}^2 \, \, ds 
    \leq \, C ( \| \vu_{0} \|^2_{H^\tau} + \| \omega_{0} \|_{H^{\sigma}}^2 ).
\end{align} 

\noindent
Applying Lemma \ref{lem5} to the inequality \eqref{fin2} with $b=1$, we find that there exists a constant $ c_2 \geq 1$ such that if
\begin{equation*}
    \| \vu_0 \|_{\dot{H}^\tau}^2 + \| \omega_{0}  \|_{H^\sigma}^2 \leq \frac 1 { c _2}
\end{equation*} 
with $ c_2  \, T_0 = 1$, then
\begin{align}
\label{cunt0}
    \| \vu_\epsilon & (t)\|^2_{H^\tau} + \| \omega_\epsilon(t) \|_{H^\sigma}^2 + \int_0^t \| \vu_\epsilon \|_{\dot{H}^{\tau+1}}^2 + \| \omega_\epsilon \|_{H^{\sigma+1}}^2 \, \, ds 
    \leq \, C ( \| \vu_{0} \|^2_{H^\tau} + \| \omega_{0} \|_{H^{\sigma}}^2 ).
\end{align} 

\noindent
That will allow us to use the following version of the Aubin--Lions theorem :

\begin{Lemma}[Aubin--Lions compactness theorem] \label{lemcom} Consider $s>0$, $q >1$ and $r <0$. Let $ (f_n)_n$ be a sequence of functions on $(0,T)\times \mathbb{R}^d$ such that, for all $T_0\in (0,T)$ and all $\varphi\in\mathcal{D}(\mathbb{R}^d)$,
  \begin{itemize}
  \item[$\bullet$]  $\varphi f_n$ is bounded in $L^2((0,T_0), H^s)$ 
  \item[$\bullet$]   $\varphi \partial_t f_n$ is bounded in $L^q((0,T_0), H^r)$ .
\end{itemize}
\noindent
Then, there exists a subsequence $(f_{n_k})$ such that for all $T_0 \in (0,T)$ and all $R_0>0$,
$$\lim_{n_k\rightarrow +\infty} \int_0^{T_0} \int_{ |  x | \leq R_0} | f_{n_k}-f_\infty |^2 \, dx\, dt=0.$$
\end{Lemma}

\noindent
For a proof of this lemma, we refer to the books \cite{BF12} and \cite{LR24}.

By \eqref{unt0} and \eqref{cunt0} we have $(\varphi \vu_\epsilon,\varphi \omega_\epsilon ) $ is bounded in $L^2 ((0,T_0),H^{\tau+1})\times L^2 ((0,T_0),H^{\sigma+1})$. 
We can verify that $(\varphi\partial_t\vu_\epsilon,\varphi\partial_t\omega_\epsilon ) $ is bounded in $ L^\alpha((0,T_0 ), H^{-s})$ for some $s \in (-\infty,0)$ and some $\alpha>1$. In fact, it is not difficult to verify that $\varphi\partial_t\vu_\epsilon$ is bounded in $L^2((0,T_0 ), H^{-3/2}$.
Moreover, from computations done to obtain the controls \eqref{pre}, we see that
$\varphi \nabla( \vv_\epsilon \otimes \omega_\epsilon )$ is uniformly bounded in $L^2((0,T_0 ), H^{ \sigma-1 } ) \subset L^2((0,T_0 ), H^{-3/2} )$ and we can deduce $\varphi\partial_t\omega_\epsilon$ is bounded in $L^2((0,T_0 ), H^{-3/2} ) $.

\noindent
Thus, by Lemma \ref{lemcom} there exists  $(\vu, \omega)$  and a  sequence $(\epsilon_k)_{k\in\mathbb{N}}$ converging to $0$ such that
\begin{equation}
\label{l2loc}
        \lim_{k\rightarrow +\infty} \int_0^{T_0} \int_{ |  y | <R}  |  \vu_{\epsilon_k}-\vu  |^2 + |  \omega_{\epsilon_k}-\omega  |^2\, dx\, ds  =0   .
\end{equation}

\noindent
Moreover, we have that $\vu_{\epsilon_k}$  converges *-weakly to $\vu $ in $L^\infty((0,T_0 ), L^2(\Phi dx))$ and $\nabla \otimes\vu_{\epsilon_k}$ converges weakly to $\nabla\otimes\vu$  in $L^2((0,T_0 ),L^2(\Phi dx))$. Then, using \eqref{l2loc} we deduce that $ \vu_{\epsilon_k}$ converges weakly to $ \vu$ in $L^3((0,T_0 ), L^3(\Phi^{\frac{3}{2}} dx))$.
Thus, we obtain that $ \vv_{\epsilon_k}\otimes \vu_{\epsilon_k}$ and $ \vv_{\epsilon_k}\otimes \omega_{\epsilon_k}$  
are weakly convergent in $(L^{6/5}L^{6/5})_{\rm loc}$ to $\vu\otimes\vu$ and $\vu\otimes\omega$ respectively, and thus in $\mathcal{D}'((0,T_0 )\times \mathbb{R}^d)$.

\noindent
By the continuity of the  Riesz transforms we find $  p_{\epsilon_k}$ is convergent to the distribution $$p= \sum_{i=1}^3\sum_{j=1}^3 R_i R_j(u_{i}u_{j}).$$

\noindent
Thus, we have obtained
\begin{equation*} \left\{ \begin{array}{ll} \vspace{2mm} 
\partial_t \vu = \Delta \vu - (\vu \cdot \nabla) \vu  - \nabla p + \frac{1}{2} \nabla \wedge \omega ,  \\ \vspace{2mm}
\partial_{t} \omega = \Delta \omega - (\vu \cdot \nabla) \omega + \frac{1}{2} \nabla \wedge \vu - \omega +\nabla ( \nabla \cdot \omega ) .
\end{array}
\right.     
\end{equation*}

Now, we want to verify the initial condition. Since $ ( \partial_t \vu, \partial_t \omega ) $  is locally in $L^1 H^{-2}$ we deduce that $(\vu, \omega)$ has representative such that  $t\mapsto (\vu, \omega)(t,.)$
is continuous from $[0,T_0 )$ to $\mathcal{D}'(\mathbb{R}^d)$ and coincides with $(\vu,\omega)(0,.)+\int_0^t \partial_t (\vu,\omega) \, ds$.

\noindent
In the sense of distributions, we can write 
\begin{align*}
(\vu,\omega)(0,.)+\int_0^t \partial_t (\vu,\omega)\, ds=&(\vu,\omega)=\lim_{k\rightarrow +\infty} (\vu_{\epsilon_k} ,\omega_{\epsilon_k} )\\
=& \lim_{k\rightarrow +\infty}
(\vu_{0,\epsilon_k},\omega_{0,\epsilon_k})+ \int_0^t \partial_t (\vu_{n_k},\omega_{n_k})\, ds \\
=&(\vu,\omega)_{0}+\int_0^t \partial_t(\vu,\omega)\, ds,
\end{align*}
hence, $(\vu,\omega)(0,.)=(\vu_{0},\omega_{0})$. Then $(\vu,\omega)$ is in fact a solution of $(M)$.

\noindent
To verify the strong convergence to the initial data, we observe that on $(0,T_0)$ we have a uniform control in $\epsilon$ and $t$, of the quantities $\| \vu_\epsilon(t,.) \|_{H^\tau}^2$ and $\| \omega_\epsilon(t,.) \|_{H^\sigma}^2$. Then, we get the estimate
$$ \|\vu_\epsilon(t,.)\|_{H^\tau}^2 + \|\omega_\epsilon(t,.)\|_{H^\sigma}^2 \leq \|\vu_{0,\epsilon}\|_{H^\tau}^2 + \|\omega_{0,\epsilon}\|_{H^\sigma}^2 + C t (1+\|\vu_0\|_{H^\tau}^{6} + \|\omega_0\|_{H^\sigma}^{6} ).$$

\noindent
We remember that $\omega_{0,\epsilon_k}$ is  strongly convergent to $\omega_0$ in $H^\sigma$. Moreover,
since $ (\vu_{\epsilon_k}, \omega_{\epsilon_k}) =
( \vu_{0, \epsilon_k}, \omega_{0, \epsilon_k} ) + \int_0^t \partial_t (\vu_{\epsilon_k}, \omega_{\epsilon_k} )  \, ds $ we have convergence of $(\vu_{\epsilon_k}, \omega_{\epsilon_k} )(t,.)$ to $(\vu, \omega)(t,.)$ in $\mathcal{D}'(\mathbb{R}^d)$, and we can deduce that it is weakly convergent in $ H^\tau \times H^\sigma $ (because of the bound in $ L^\infty H^\tau \times L^\infty H^\sigma$), so that
$$ \|\vu(t,.)\|_{H^\tau}^2 + \|\omega(t,.)\|_{H^\sigma}^2 \leq \|\vu_{0}\|_{H^\tau}^2 + \|\omega_{0}\|_{H^\sigma}^2  + C t (1 + \|\vu_0\|_{H^\tau}^{6} + \|\omega_0\|_{H^\sigma}^{6} ) . $$

\noindent
Taking the limit when $t$ goes to $0$, this remark implies
$$\limsup_{t\rightarrow 0}  \|\vu(t,. )\|_{H^\tau}^2 + \|\omega(t,. )\|_{H^\sigma}^2 \leq  \|\vu_{0} \|_{H^\tau}^2 + \|\omega_{0} \|_{H^\sigma}^2 .$$

\noindent
For the reciprocal inequality, we recall that $\vu$ is weakly continuous in $H^\tau$ and $\omega$ is weakly continuous in $H^\sigma$, therefore
\begin{equation*}
    \|\vu_{0} \|_{H^\tau}^2 \leq   \liminf_{t\rightarrow 0}  \|\vu(t,. )\|_{H^\tau}^2 \phantom{spa} \text{and} \phantom{spa} \|\omega_{0} \|_{H^\sigma}^2 \leq  \liminf_{t\rightarrow 0}  \|\omega(t,. )\|_{H^\sigma}^2 .
\end{equation*}
  
\noindent
Thus, we can conclude that $$ \sqrt{\|\vu_{0} \|_{H^\tau}^2 + \|\omega_{0} \|_{H^\sigma}^2} =   \lim_{t\rightarrow 0} \sqrt{ \|\vu(t,. )\|_{H^\tau}^2 + \|\omega(t,. )\|_{H^\sigma}^2 }. $$

\noindent
It allows us to turn the weak convergence into a strong convergence in the Hilbert space $H^\tau \times H^\sigma$. It finish the proof.

\section*{Conflict of Interest}
The author declare that he have no conflict of interest.

\section*{Data Availability Statement}
Data sharing is not applicable to this article as no datasets were generated or analyzed during the current study.

\section*{Acknowledgment}
PF is supported by the Basque Government through the BERC 2022-2025 program and by the Spanish State Research Agency through BCAM Severo Ochoa CEX2021-001142

\end{document}